\documentclass[11pt,leqno]{amsart}
\usepackage{amssymb,verbatim,enumerate,ifthen}
\usepackage[mathscr]{eucal}
\usepackage[utf8]{inputenc}
\usepackage[T1]{fontenc}
\def\H{\mathbb{H}}
\def\N{\mathbb{N}}
\def\R{\mathbb{R}}

\def\H{\mathscr{H}}

\newtheorem{theorem}{Theorem}[section]
\newtheorem*{theorem*}{Theorem}
\def\Thm#1#2{\ifthenelse{\equal{#1}{*}}{\begin{theorem*}#2\end{theorem*}}
             {\begin{theorem}\label{T#1}#2\end{theorem}}}
\newtheorem{Atheorem}{Theorem}

\def\thm#1{Theorem~\ref{T#1}}
\newtheorem{proposition}[theorem]{Proposition}
\newtheorem*{proposition*}{Proposition}
\def\Prp#1#2{\ifthenelse{\equal{#1}{*}}{\begin{proposition*}#2\end{proposition*}}
{\begin{proposition}\label{P#1}#2\end{proposition}}}
\def\prp#1{Proposition~\ref{P#1}}

\newtheorem{corollary}[theorem]{Corollary}
\newtheorem*{corollary*}{Corollary}
\def\Cor#1#2{\ifthenelse{\equal{#1}{*}}{\begin{corollary*}#2\end{corollary*}}
             {\begin{corollary}\label{C#1}#2\end{corollary}}}

\newtheorem{lemma}[theorem]{Lemma}
\newtheorem*{lemma*}{Lemma}
\def\Lem#1#2{\ifthenelse{\equal{#1}{*}}{\begin{lemma*}#2\end{lemma*}}
             {\begin{lemma}\label{L#1}#2\end{lemma}}}
\def\lem#1{Lemma~\ref{L#1}}

\theoremstyle{definition}
\newtheorem{remark}[theorem]{Remark}
\newtheorem*{remark*}{Remark}
\def\Rem#1#2{\ifthenelse{\equal{#1}{*}}{\begin{remark}\rm #2\end{remark}}
             {\begin{remark}\label{R#1}\rm #2\end{remark}}}

\newtheorem{example}[theorem]{Example}
\newtheorem*{example*}{Example}
\def\Exa#1#2{\ifthenelse{\equal{#1}{*}}{\begin{example*}\rm #2\end{example*}}
             {\begin{example}\label{Ex#1}\rm #2\end{example}}}

\def\eq#1{{\rm(\ref{E#1})}}
\def\Eq#1#2{\ifthenelse{\equal{#1}{*}}
  {\begin{equation*}\begin{aligned}#2\end{aligned}\end{equation*}}
  {\begin{equation}\begin{aligned}\label{E#1}#2\end{aligned}\end{equation}}}

\begin{document}
\begin{flushright}
\end{flushright}
\vspace{5mm}

\date{\today}

\title{Convex Sequence and Convex Polygon}

\author[A. R. Goswami]{Angshuman R. Goswami}
\address[A. R. Goswami]{Department of Mathematics, University of Pannonia, 
H-8200 Veszprem, Hungary}
\email{goswami.angshuman.robin@mik.uni-pannon.hu}
\author[István Szalkai]{István Szalkai}
\address[Szalkai István]{Department of Mathematics, University of Pannonia, 
H-8200 Veszprem, Hungary}
\email{szalkai.istvan@mik.uni-pannon.hu}
\subjclass[2000]{Primary: 26A48; Secondary: 26A12, 26A16, 26A45}
\keywords{convex sequence; convex function; convex polygon}

\begin{abstract}
In this paper, we deal with the question; under what conditions the points 
\newline
$P_i(xi,yi)$ $(i = 1,\cdots, n)$ form
a convex polygon provided $x_1 < \cdots < x_n$ holds. One of the main findings of the paper can be stated as follows\\

"Let $P_1(x_1,y_1),\cdots ,P_n(x_n,y_n)$ are $n$ distinct points ($n\geq3$) with $x_1<\cdots<x_n$. Then $\overline{P_1P_2},\cdots \overline{P_nP_1}$ form a convex $n$-gon that lies in the half-space
\Eq{*}{
\underline{\H}=\bigg\{(x,y)\big|\quad x\in\R \quad \mbox{and} \quad  y\leq y_1+\bigg(\dfrac{x-x_1}{x_n-x_1}\bigg)(y_n-y_1)\bigg\}\subseteq{\R^{2}}
}
if and only if the following inequality holds
\Eq{*}{
\dfrac{y_i-y_{i-1}}{x_i-x_{i-1}} \leq \dfrac{y_{i+1}-y_{i}}{x_{i+1}-x_{i}} \quad \quad \mbox{for all} \quad \quad i\in\{2,\cdots,n-1\}
."}
Based on this result, we establish a linkage between the property of sequential convexity and convex polygon. We show that in a plane if any $n$ points are scattered in such a way that their horizontal and vertical distances preserve some specific monotonic properties; then those points form a $2$-dimensional convex polytope. \\

Various definitions, backgrounds, motivations, findings, and other important matters are discussed in the introduction section.
\end{abstract}

\maketitle

\section*{Introduction} 
Throughout this paper $\N$, $\R$, and $\R_+$ denote the set of natural, real, and positive numbers respectively. $\R^2$ is used to indicate the usual $2$-dimensional plane.\\

A sequence $\big<u_i\big>_{i=0}^{\infty}$ is said to be convex if it satisfies the following inequality
\Eq{*}{
2u_i\leq u_{i-1}+u_{i+1}\quad\mbox{for all} \quad i\in\N.
}  
In other words, $\big<u_i\big>_{i=0}^{\infty}$ possesses sequential convexity if the sequence $\big<u_i-u_{i-1}\big>_{i=1}^{\infty}$ is increasing. If the converse of the above inequality holds, we will term $\big<u_i\big>_{i=0}^{\infty}$ as a concave sequence. \\

The very first evidence where the terminology of sequential convexity was used is the book of
\cite{Mitrinovicc}. Since then many important results have been discovered in this direction; such as establishing a discrete version of Hermite-Hadamard type inequality, Ulam's type stability theorems, applications in the field of trigonometric functions, generalization of sequential convexity to the higher order and in an approximate sense. The details of these can be found in the papers [\cite{Mitrinovicc}-\cite{Krasniqi}] and the mentioned references there.\\

In discrete geometry, there are many results that primarily mention about the scattered random points and the underlying convex geometry. One of the familiar examples of it is well-known Radon's theorem. Helly's and Caratheodory's results also indirectly deal with the same. The background, origin, generalization, and other research development related details can be found in the papers \cite{Radon, Helly, Carathéodory, Ernst, Imre} and in the book \cite{Barvinok}. There are many tempting open problems that deal with questions regarding the possibilities of existing a specifically shaped convex body in a higher dimensional space provided the scattered point bounds with some specific patterns or numbers. For instance, the famous Erd\"os–Szekeres conjecture is still unsolved even after almost 90 years of first mention. To better understand the problem, we can look into \cite{Erdos}. \\

The purpose of this paper is to investigate a necessary and sufficient condition that $n$ distinct points in $\R^2$ namely $P_1(x_1,y_1),\cdots ,P_n(x_n,y_n)$  with $x_1<\cdots<x_n$ must follow in order to form a convex $n-$gon. It turns out that if the following inequality holds 
\Eq{*}{
\dfrac{y_i-y_{i-1}}{x_i-x_{i-1}} \leq \dfrac{y_{i+1}-y_{i}}{x_{i+1}-x_{i}} \quad \quad \mbox{for all} \quad \quad i\in\{2,\cdots,n-1\}}
then the line segments $\overline{P_1P_2},\cdots \overline{P_nP_1}$ form a convex polygon with $n$-sides that lies below $\overline{P_nP_1}$. The converse implication of this statement can also be established. Similarly under the same assumptions the reverse inequality holds if and only if the the convex $n$-gon lies above the line segment $\overline{P_nP_1}$.\\

Based on this result, we derived some other interesting findings. We assume some sequential convexity properties on $\big<x_i\big>_{i=1}^{n}$ and $\big<y_i\big>_{i=1}^{n}$ as follows
\begin{enumerate}[(i)]
\item $\big<x_i\big>_{i=1}^{n}$ is strictly increasing and concave;\\
\item $\big<y_i\big>_{i=1}^{n}$ is increasing and convex.
\end{enumerate}
Then $P_1,\cdots,P_n$ are the vertices of a $n$-convex polygon.  A more generalized version of this result is also presented.\\

Having prior knowledge of the vertices of a convex polygon often reduces lot of mathematical and computational tasks. In linear programming problems, the optimized value of the cost function always lies in one of the vertices of the constrains formulated convex polygon. In computational geometry, efficient algorithms are highly dependent upon the extreme points of convex $n$-gons. In computer graphics, various rotation, translation and orientation related techniques are performed at the end points of a convex polygon. \\

In the way of proving our results, we establish several lemmas, and propositions related to fractional inequality, convex sequence, and convex function theory.
\section{Main Results}
Our first result shows a very important fractional inequality. Later this inequality is going to be used extensively to establish some of the results. This inequality is also mentioned in one of our recently submitted papers. However, for readability purpose, we mentioned the statement along with the proof.
\Lem{22}{Let $n\in\N$ be arbitrary. Then for any  $a_1,\cdots,a_n\in\R$ and $b_1,\cdots,b_n\in\R_+$; the following inequalities hold
\Eq{111}{
\min\bigg\{\dfrac{a_1}{b_1},\cdots, \dfrac{a_n}{b_n}\bigg\}
\leq\dfrac{a_1+\cdots+a_n}{b_1+\cdots+b_n}\leq \max\bigg\{\dfrac{a_1}{b_1},\cdots, \dfrac{a_n}{b_n}\bigg\}.}}
\begin{proof}
We will prove the theorem by using mathematical induction. For $n=1$, there is nothing to prove. For  $a_1,a_2\in\R$ and $b_1,b_2\in\R_+$, without loss of generality, we can assume that $\dfrac{a_1}{b_1}\leq \dfrac{a_2}{b_2}$; which is equivalent to
\Eq{*}{
a_1(b_1+b_2)\leq (a_1+a_2) b_1 \quad\mbox{and}\quad (a_1+a_2)b_2\leq (b_1+ b_2)a_2.
}
The two inequalities above together yield the following 
\Eq{67}
{\dfrac{a_1}{b_1}\leq \dfrac{a_1+a_2}{b_1+b_2}\leq \dfrac{a_2}{b_2}}
which validates \eq{111} for $n=2.$ Now we assume that the statement is true for a $n\in\N.$
 
Let $a_1,\cdots,a_n,a_{n+1}\in\R$ and $b_1,\cdots,b_n,b_{n+1}\in\R_+$. Utilizing \eq{67} and our induction assumption \eq{111}, we can compute the following two inequalities 
\Eq{*}{
&\min\bigg\{ \dfrac{a_1}{b_1}, \cdots ,\dfrac{a_n}{b_n}, \dfrac{a_{n+1}}{b_{n+1}} \bigg\}\leq \min \bigg\{ \dfrac{a_1 +\cdots+a_n}{b_1+\cdots+b_n},\dfrac{a_{n+1}}{a_{n+1}}\bigg\}\leq \dfrac{a_1+\cdots+a_n+a_{n+1}}{b_1+\cdots+b_n+b_{n+1}}
\\
&\qquad\qquad\qquad\qquad\qquad\qquad\qquad\qquad\qquad\mbox{and}\\
&\dfrac{a_1+\cdots+a_n+a_{n+1}}{b_1+\cdots+b_n+b_{n+1}}\leq\max\bigg\{ \dfrac{a_1+\cdots+a_n}{b_1+\cdots+b_n},\dfrac{a_{n+1}}{a_{n+1}}\bigg\}\leq \max\bigg\{ \dfrac{a_1}{b_1}, \cdots ,\dfrac{a_n}{b_n}, \dfrac{a_{n+1}}{b_{n+1}} \bigg\}.}
This justifies \eq{111} for any $n\in\N$ and completes our proof.
\end{proof}
We can obtain several mean inequalities as a consequence of the above result. 
\newline If $b_1=\cdots=b_n=1$; then \eq{111} turns into the  standard arithmetic mean inequality which is represented by as follows
\Eq{*}{
\min\{a_1,\cdots,a_n\}\leq \dfrac{a_1+\cdots+a_n}{n}\leq \max\{a_1,\cdots,a_n\}.
}
On the other hand, we can consider $b_i=1/x_i$ and $a_i=1$ for all $i\in\{1,\cdots,n\}$. Upon substituting these in \eq{111}, we get the Harmonic mean inequality for the positive numbers $x_1,\cdots,x_n$ that can be formulated as
\Eq{*}{
\min \{x_1,\cdots, x_n\}\leq \dfrac{n}{\dfrac{1}{x_1}+ \cdots + \dfrac{1}{x_n}}\leq \max \{x_1, \cdots, x_n \}.
} 
Before moving to the main result of this section, there are some notions and terminologies that we need to recall.
There are several ways to represent a convex function. Besides the standard definition of convexity; for any given function, the monotonic property of the associated slope function can also be used to determine convexity. In other words,
A function $f:I\to\R$ is said to be convex if for any $x',$ $x$ and $x"\in I$ with $x'<x<x''$, the below mentioned inequality holds
\Eq{20011}{ 
\dfrac{f(x)-f(x')}{x-x'}\leq \dfrac{f(x")-f(x)}{x"-x}.
}
The other concept we are going to utilize is the epigraph of a function. For a function $f:I\to\R$; the notion of epigraph can be formulated as follows
\Eq{*}{
epi(f)=\Big\{(x,y):\, \, f(x)\leq y,\,\, x\in I \Big\}.
}
One of the basic characterizations of a convex function can be stated as "A function $f$ is convex if and only if  $epi(f)$ is a convex set."\\

Now we have all the required tools to proceed to state the first theorem.
\Thm{991}{Let $P_1(x_1,y_1),\cdots ,P_n(x_n,y_n)$ are $n$ points($n\geq3$) with $x_1<\cdots<x_n$. Then $\overline{P_1P_2},\cdots \overline{P_nP_1}$ form a convex $n$-gon in the below mentioned half-space
\Eq{12}{
\underline{\H}=\bigg\{(x,y)\big|\quad x\in\R \quad \mbox{and} \quad  y\leq y_1+\bigg(\dfrac{x-x_1}{x_n-x_1}\bigg)(y_n-y_1)\bigg\}\subseteq{\R^{2}}
}
if and only if the following inequality holds 
\Eq{13}{
\dfrac{y_i-y_{i-1}}{x_i-x_{i-1}} \leq \dfrac{y_{i+1}-y_{i}}{x_{i+1}-x_{i}} \quad \quad \mbox{for all} \quad \quad i\in\{2,\cdots,n-1\}
.}
}
\begin{proof}
There are several steps involve in the proof.
To prove the theorem, we assume that \eq{13} is valid. We define the function $f:[x_1,x_n]\to\R$ as follows
\Eq{7878}{f(x):=ty_{i}+(1-t)y_{i+1}\quad \mbox{where} \quad x:=tx_i+(1-t)x_{i+1} \quad (t\in [0,1]\, \mbox{  and  }\, i\in \{1,\cdots,n-1\}.}
From the construction, it is clearly visible that $f$ is formulated by joining total $n-1$ consecutive line segments that are defined in between the points $x_{i}$ and $x_{i+1}$ for all $i\in \{1,\cdots,n-1\}$. we term these as functional line segments of $f$. For the proof, first we are going to establish that the function $f$ is convex. But before that we need to validate the statement below
\Eq{200}{
&\mbox{\textbf{Slope of the functional line segment(s) of $f$ that contains $(x',f(x'))$}}\\
&\leq \mbox{\textbf{Slope of the line joining the points $(x',f(x'))$ and $(x'',f(x''))$}}\\
&\leq \mbox{\textbf{Slope of the functional line segment(s) of $f$ that contains $(x'',f(x''))$}.}
}
If both the points $(x',f(x'))$ and $(x'',f(x''))$ lie in the same functional line segment, the statement is obvious.  

Next we consider the case when $x'\in[x_{i-1},$ $x_i[$ and $x''\in]x_{i},x_{i+1}]$, ($i\in\{2,\cdots,n-1\}$) that is $x'$ and $x''$; lie in two consecutive intervals.
Utilizing \eq{13} and basic geometry of slopes in straight lines, we can compute the inequality below
\Eq{2001}{
\dfrac{f(x_i)-f(x')}{x_i-x'}=\dfrac{f(x_i)-f(x_{i-1})}{x_i-x{i-1}}\leq \dfrac{f(x_{i+1})-f(x_i)}{x_{i+1}-x_i}=\dfrac{f(x")-f(x_i)}{x"-x_i}
.} 
The expression $\dfrac{f(x")-f(x')}{x"-x'}$ can also be written as $\dfrac{f(x")-f(x_i)+f(x_i)-f(x')}{(x"-x_i)+(x_i-x')}$. 
This along with \eq{67} and \eq{2001} yields
\Eq{*}{
\dfrac{f(x_i)-f(x')}{x_i-x'}\leq \dfrac{f(x")-f(x')}{x"-x'}\leq\dfrac{f(x")-f(x_i)}{x"-x_i}
.} 
and validates the statement \eq{200} for this particular case.\\
Finally we assume $x'\in[x_j,x_{j+1}]$, $x"\in[x_k,x_{k+1}]$ where $j\in\{1,\cdots,n-2\}$ and $k\in\{3,\cdots,n\}$ such that $k-j\geq 2$.
First using \eq{111} of \lem{22} and then by applying \eq{13}; we can obtain the inequality below
\Eq{*}{
\dfrac{f(x")-f(x')}{x"-x'}&= \dfrac{\Big(f(x")-f(x_k)\Big)+\Big(f(x_k)-f(x_{k-1})\Big)+\cdots +\Big(f(x_{j+1})-f(x')\Big)}{\Big(x"-x_k\Big)+\Big(x_k-x_{k-1}\Big)+\cdots +\Big(x_{j+1}-x'\Big)}\\
&\leq\max\bigg\{\dfrac{f(x")-f(x_{k})}{x"-x_k},\dfrac{f(x_k)-f(x_{k-1})}{x_k-x_{k-1}},\cdots,\dfrac{f(x_{j+1})-f(x_j)}{x_{j+1}-x_j}\bigg\}\\
&=\dfrac{f(x")-f(x_{k})}{x"-x_k}=\dfrac{f(x_{k+1})-f(x_{k})}{x_{k+1}-x_k}.}
Similarly, we can compute the following inequality as well
\Eq{*}{
\dfrac{f(x_{j+1})-f(x_{j})}{x_{j+1}-x_j}&=\dfrac{f(x_{j+1})-f(x')}{x_{j+1}-x'}\\
&=\min\bigg\{ \dfrac{f(x_{j+1})-f(x')}{x_{j+1}-x'},\dfrac{f(x_{j+2})-f(x_{j+1})}{x_{j+2}-x_{j+1}},\cdots,\dfrac{f(x")-f(x_k)}{x"-x_k}\bigg\}\\
&\leq \dfrac{\Big( f(x_{j+1})-f(x')\Big)+\Big(f(x_{j+2})-f(x_{j+1})\Big)+\cdots+\Big(f(x")-f(x_i)\Big)}{\Big(x_{j+1}-x'\Big)+\Big(x_{j+2}-x_{j+1})\Big)+\cdots+\Big(x'"-x_k \Big)} 
\\
&=\dfrac{f(x")-f(x')}{x"-x'}.}
The above two inequalities establishes the statement \eq{200}.\\
We are now ready to show that $f$ is convex. we assume  $x',x,x"\in[x_1,x_n]$ with $x'<x<x"$. More specifically, $x'\in[x_j,x_{j+1}]$,$x\in[x_i,x_{i+1}]$ and $x"\in[x_k,x_{k+1}]$ for fixed $i,j$ and $k\in\{1,\cdots,n-1\}$.  Then first using the last part of inequality in \eq{200} and then applying the initial part of the same inequality, we obtain the following
\Eq{*}{
\dfrac{f(x)-f(x')}{x-x'}&\leq\dfrac{f(x_{i+1})-f(x_i)}{x_{i+1}-x_{i}}=\dfrac{y_{i+1}-y_i}{x_{i+1}-x_{i}}\\
&\quad\quad\quad\quad\,\,\mbox{and}\\
\dfrac{y_{i+1}-y_i}{x_{i+1}-x_{i}}&=\dfrac{f(x_{i+1})-f(x_i)}{x_{i+1}-x_{i}}\leq \dfrac{f(x")-f(x)}{x"-x}.
}
Combining the above two inequalities, we arrive at \eq{20011}. It shows that the function $f$ possesses convexity and also implies $epi(f)$ is an unbounded convex polygon. Due to the convexity property of $f$, the line segment $\overline{P_nP_1}$ that forms by joining $(x_1,y_1)$ and $(x_n,y_n)$  entirely lies in $epi(f)$. 
The extension of $\overline{P_nP_1}$ creates two closed convex 
half-spaces. One of these is defined as in \eq{12}.

All these yields, $\underline{\H}\cap epi(f)$ gives the convex polygon $\overline{P_1P_2}$, $\overline{P_2P_3},\cdots ,\overline{P_nP_1}$ that lies in $\underline{\H}$.\\

Conversely, suppose that under the strict monotonic assumptions on $x_i's$, the $n$ distinct points $P_1,\cdots,P_n$ form a convex polygon. A convex polyhedron is nothing but the intersection of a finite number of hyper-planes. By neglecting the underlying hyperplane due to the extension of $\overline{P_nP_1}$; we will end up in an unbounded convex polygon formed by the line segments $\overline{P_1P_2}, \cdots ,\overline{P_{n-1}P_n}.$ In other words; this unbounded convex set is just the epigraph of the function $f:[x_1,x_n]\to\R$; that was defined in \eq{7878}. Convexity of $epi(f)$ implies $f$ is a convex function. Since, $P_1(x_1,y_1),\cdots ,P_n(x_n,y_n)\in gr(f)$; by \eq{20011} we can also conclude that \eq{13} holds.
This completes the proof of the statement.
\end{proof}
\begin{remark}
In the above theorem, the strict monotonicity of 
the sequence $\big<x_i\big>_{i=1}^{n}$ can be relaxed for the end points.
For instance, instead of strict increasingness; we can assume that the elements of the sequence
$\big<x_i\big>_{i=1}^{n}$ satisfies the following inequality
$$ x_1\leq x_2<x_3<\cdots<x_{n-1}\leq x_n.$$ In the case of $x_1=x_2$, the points $P_1,\cdots,P_n$ still configure a $n$ sided convex polygon provided \eq{13} holds for all $i\in\{2,\cdots,n-1\}$. Similar to the proof of \thm{991}; one can show that, by joining of $n-1$ distinct points $(x_2,y_2),\cdots,(x_n,y_n) $ we can formulate the convex function $f$. Then $\overline{\H}\cap epi(f)$ results in a convex set. This set is bounded from below by the $n-2$ line segments that form the function $f$. The line segment joining the points $(x_1,y_1)$ and $(x_n,y_n)$ gives the upper bound. And finally, the line segment that passes through $(x_1,y_1)$ and $(x_2,y_2)$  also lies in the set. In other words, we ended up in a convex polygon which is formed by the line segments 
$\overline{P_1P_2}, \cdots ,\overline{P_{n}P_1}.$ Similar conclusion can be drawn by considering equality in the right most or simultaneously in both end points of the sequence $\big<x_i\big>_{i=1}^{n}.$ 
\end{remark}
We can also establish the following theorem. The proof of it is analogous to \thm{991}. Hence, the proof is not included.
\Thm{992}{Let $P_1(x_1,y_1),\cdots ,P_n(x_n,y_n)$ are $n$ points($n\geq3$) with $x_1<\cdots<x_n$. Then $\overline{P_1P_2},\cdots \overline{P_nP_1}$ form a convex $n$-gon in the below mentioned half-space
\Eq{*}{
\underline{\H}=\bigg\{(x,y)\big|\quad x\in\R \quad \mbox{and} \quad   y_1+\bigg(\dfrac{y_n-y_1}{x_n-x_1}\bigg)(x-x_1)\leq y\bigg\}\subseteq{\R^{2}}
}
if and only if the following inequality holds 
\Eq{*}{\dfrac{y_{i}-y_{i-1}}{x_{i}-x_{i-1}}\geq \dfrac{y_{i+1}-y_{i}}{x_{i+1}-x_{i}} \quad \quad \mbox{for all} \quad \quad i\in\{2,\cdots,n-1\}
.}
}
In the next section, we are going to see how sequential convexity is linked with a convex polygon.
\section{The subsequent results}
The results of this section heavily depend upon the previous section's findings.
The initial proposition resembles the slope property of a convex function.
\Prp{1}{Suppose $\big<x_i\big>_{i=0}^{\infty}$ is a strictly increasing and concave sequence.  While $\big<y_i\big>_{i=1}^{\infty}$ is a increasing and convex sequence. Then for any $i\in\N$, the following discrete functional inequality holds
\Eq{1}{
\dfrac{y_i-y_{i-1}}{x_i-x_{i-1}}\leq \dfrac{y_{i+1}-y_{i}}{x_{i+1}-x_{i}}.}
}
\begin{proof}
By our assumptions $0<x_{i+1}-x_i\leq x_i-x_{i-1}$ and $0\leq y_{i}-y_{i-1}\leq y_{i+1}-y_{i}$. Multiplying these two inequalities side by side and rearranging the terms of the resultant inequality, we obtain \eq{1}.
\end{proof}
In the above proposition, the respective positive/non-negative conditions in sequences $\big<x_i-x_{i-1}\big>_{i=1}^{\infty}$ and $\big<y_i-y_{i-1}\big>_{i=1}^{\infty}$ cannot be compromised. One can easily verify this fact by relaxing this crucial condition. Similar to the above proposition, we can also show the next result.
\Prp{2}{Suppose $\big<x_i\big>_{i=0}^{\infty}$ is a strictly increasing convex sequence. While $\big<y_i\big>_{i=1}^{\infty}$ is a decreasing  convex sequence.  Then for any $i\in\N$, the discrete functional inequality \eq{1} is satisfied.}

Now we can propose the first theorem of this section. The establishment of it is similar to the first part of \thm{991}; hence only the statement is mentioned. 

\Thm{15}{Let $P_1,\cdots,P_n$ are $n$ points with the respective coordinates $(x_1,y_1),\cdots,(x_n,y_n)$  scattered in the standard $\R^2$  such that the sequences $\big<x_i\big>_{i=1}^{n}$ is strictly monotone and concave; while  $\big<y_i\big>_{i=1}^{n}$ is convex and increasing . Then $\overline{P_1P_2}, \cdots ,\overline{P_{n-1}P_n}.$  form a convex polygon. }
\begin{proof}
This theorem is a direct consequence of \prp{1} and \thm{991}.
\end{proof}
The converse of the above theorem is not necessarily always true. We consider the three vertices of the $ {\displaystyle \bigtriangleup }P_1P_2P_3$ as $P_1(0,0),P_2(2,2)$ and $P_3(3,1)$. One can easily observe that both the $X$ and $Y$-coordinated sequences are strictly concave which shows reverse implication is not valid. \\

The next result is similar to the above one. By using \prp{2} and  \thm{991}; we can establish it. 
\Thm{16}{Let $P_1,\cdots,P_n$ are $n$ points with the respective coordinates $(x_1,y_1),\cdots,(x_n,y_n)$  scattered in the standard $\R^2$  such that the sequences $\big<x_i\big>_{i=1}^{n}$ is strictly monotone and convex; while  $\big<y_i\big>_{i=1}^{n}$ is convex with $y_{i+1}-y_{i}\leq 0$ for all $i\in\{1,\cdots,n-1\}$. Then $P_1,\cdots P_n$ form a convex polygon.}

A more general result can be formulated by combining these two theorems. But Before proceeding, we must go through a proposition.
\Prp{3}{Let $\big<u_i\big>_{i=1}^{n}$ be a convex sequence. Then there exists an element $m\in\N\cap[1,n]$, that satisfies at least one of the following 
\Eq{60}{
u_i\leq u_m \quad\mbox{for all $i<m$}\quad \quad\quad\mbox{or}\quad \quad\quad 
u_m\leq u_i \quad\mbox{for all $m<i$}
.}}
\begin{proof}
The sequential convexity of $\big<u_i\big>_{i=1}^{n}$ implies increasingness of the sequence $\big<u_{i+1}-u_{i}\big>_{i=1}^{n-1}.$
If all the terms in this monotone sequence are either non-negative or non-positive, then one can easily validates \eq{60}. If not, then there exists a $m\in\N\cap]1,n[$ such that the following inequalities holds 
\Eq{*}{
u_{m+1}-u_{m}\geq 0 \quad\quad \quad\quad\mbox{or}\quad \quad\quad 
u_{m}- u_{m-1}\leq 0.
}
This together with non-decreasingness property of the sequence $\big<u_{i+1}-u_{i-1}\big>_{i=1}^{n-1}$ yields \eq{60} and completes the proof of the statement.
\end{proof}

The next theorem generalizes our previous results  \thm{15} and \thm{16}. Hence just a scratch of the proof is mentioned.
\Thm{17}{Let $P_1,\cdots,P_n$ are $n$ distinct points with the respective coordinates $(x_1,y_1),\cdots,(x_n,y_n)$  scattered in the standard $\R^2$.  The sequences $\big<y_i\big>_{i=1}^{n}$ is convex such that 
$\underset{1\leq i\leq n}\min y_i=y_m$. Suppose the sequence $\big<x_i\big>_{i=1}^{n}$ is increasing and the sub-sequences $\big<x_i\big>_{i=1}^{m}$ and $\big<x_i\big>_{i=m}^{n}$ are sequentially convex and concave  respectively. Then the points $P_1,\cdots P_n$ form a convex polygon. }
\begin{proof}
We construct the function $f$ as in \eq{7878}. From there by utilizing \lem{22}, \prp{1}, \prp{2}, \thm{15} and \thm{16} we can conclude \eq{200}. Analogous to \thm{991}; it leads us to the establishment that the function $f$ is convex. Finally from the $epi(f)$; we can get the desired result. 
\end{proof}
Off course, as mentioned in the remark of \thm{15}; one can discuss relaxing strict monotonic property of the sequence $\big<x_i\big>_{i=1}^{n}$ in its extreme points to simply increasingness.\\

This investigation also raises several interesting new problems and challenges. For instance, one obvious task is generalizing this concept to any finite dimension. It leads us to question, in higher dimensions, what must be the coordinate-oriented inequalities need to be satisfied by the scattered points in order to obtain a convex polytope?

Another area of discussion is the possible outcome in $\R^{2}$; if the coordinates of the points follow higher order sequential convexity/concavity properties. Are we still able to extract an underlying convex polygon out of it or a complex interesting geometric figure?

\bibliographystyle{plain}

\end{document}